\documentclass [11pt]{article}
\usepackage{amsthm}
\usepackage{amsmath}
\usepackage{amssymb}
\usepackage{enumerate}
\usepackage{epsfig}

\newcounter{contador}
\newtheorem{propo}[contador]{Proposition}
\newtheorem{teo}[contador]{Theorem}
\newtheorem{lem}[contador]{Lemma}

\newtheorem{example}[contador]{Example}

\newcommand{\rec}{\noindent} 
\newcommand{\dem}{\rec {\it Proof. }} 
\renewcommand{\qed}{\, \hfill\rule[-1mm]{2mm}{3.2mm}} 

\newcommand{\Qm}{{Q^+}}

\newcommand{\sign}{{\rm sign}}

\newcommand{\enya}{${\rm \tilde{n}}$}

\newcommand{\R}{{\mathbb R}}

\newcommand{\U}{{\cal{U}}}

\title{On two and three periodic Lyness difference equations\footnote{{\bf Acknowledgements}.
GSD-UAB and CoDALab Groups are supported by the Government of Catalonia through
the SGR program. They are also supported by MCYT through grants MTM2008-03437
(first and second authors) and DPI2008-06699-C02-02 (third author).}}

\author{Anna Cima$^{(1)}$, Armengol Gasull$^{(1)}$ and V\'{\i}ctor Ma\~{n}osa $^{(2)}$
\\*[.1truecm]
{\small \textsl{$^{(1)}$ Dept. de Matem\`{a}tiques, Facultat de Ci\`{e}ncies,}}
\\*[-.25truecm] {\small \textsl{Universitat Aut\`{o}noma de Barcelona,}}
\\*[-.25truecm] {\small \textsl{08193 Bellaterra, Barcelona, Spain}}
\\*[-.25truecm] {\small \textsl{cima@mat.uab.cat, gasull@mat.uab.cat}}
\\*[-.25truecm]
\\*[-.25truecm] {\small \textsl{$^{(2)}$ Dept. de Matem\`{a}tica Aplicada III (MA3),}}
\\*[-.25truecm] {\small \textsl{Control, Dynamics and Applications Group (CoDALab)}}
\\*[-.25truecm] {\small \textsl{Universitat Polit\`{e}cnica de Catalunya (UPC)}}
\\*[-.25truecm] {\small \textsl{Colom 1, 08222 Terrassa, Spain}}
\\*[-.25truecm] {\small \textsl{victor.manosa@upc.edu}}}
\date{}

\begin{document}
\maketitle


\begin{abstract}
We describe the sequences $\{x_n\}_n$ given by the  non-autonomous
second order Lyness difference equations
$x_{n+2}=(a_n+x_{n+1})/x_n,$ where $\{a_n\}_n$ is either a
2-periodic or a 3-periodic sequence of positive values and the
initial conditions $x_1,x_2$ are as well positive. We also show an
interesting phenomenon of the discrete dynamical systems associated
to some of these difference equations: the existence of one
oscillation of their associated rotation number functions. This
behavior does not appear for the autonomous Lyness difference
equations.
\end{abstract}

\rec {\sl 2000 Mathematics Subject Classification:} \texttt{39A20, 39A11}

\rec {\sl Keywords:} Difference equations with periodic coefficients, circle
maps, rotation number.\newline

\section{Introduction and main result}

This paper fully describes the sequences given by the non-autonomous
second order Lyness difference equations
\begin{equation}\label{eq}
x_{n+2}\,=\,\frac{a_n+x_{n+1}}{x_n},
\end{equation}
where $\{a_n\}_n$ is a $k$-periodic sequence taking positive values,
$k=2,3,$ and the initial conditions $x_1,x_2$ are as well positive.
This question is proposed in \cite[Sec. 5.43]{CL}. Recall that
non-autonomous recurrences appear for instance as population models
with a variable structure affected by some seasonality
\cite{ES1,ES2}, where $k$ is the number of seasons. Some dynamical
issues of similar type of equations have been studied in several
recent papers \cite{BHS,CJK,deA,JKN,KN,FJL,GKL}.

Recall that when $k=1,$ that is $a_n=a>0,$ for all $n\in\mathbb{N}$,
then (\ref{eq}) is the famous Lyness recurrence which is well
understood, see for instance \cite{BR,Z}. The cases $k=2,3$ have
been already studied and some partial results are established. For
both cases it is known that the solutions are persistent near a
given $k$-periodic solution, which is stable. This is proved by
using some known invariants, see \cite{KN, FJL, GKL}. Recall that in
our context it is said that a solution $\{x_n\}_n$ is persistent if
there exist two real positive constants $c$ and $C,$ which depend on
the initial conditions, such that for all $n\ge1, 0<c<x_n<C<\infty.$
We prove:

\begin{teo}\label{main} Let $\{x_n\}_n$ be any  sequence defined by
(\ref{eq}) and $k\in\{2,3\}$. Then it is persistent. Furthermore, either
\begin{enumerate}

\item [(a)] the sequence $\{x_n\}_n$ is periodic, with period a multiple of $k$; or

\item [(b)] the sequence $\{x_n\}_n$
densely fills one or two (resp. one, two or three) disjoint intervals of
$\mathbb{R}^+$ when $\{a_n\}_n$ is 2-periodic (resp. 3-periodic). Moreover it
is possible by algebraic tools to distinguish which is the situation.
\end{enumerate}

\end{teo}

Our approach to describe the sequences $\{x_n\}_n$ is based on the
study of the natural dynamical system associated to (\ref{eq}) and
on the results of \cite{CGM}. The main tool that allows to
distinguish the number of intervals for the adherence of the
sequences $\{x_n\}_n$ is the computation of several resultants, see
Section~\ref{provaA}.

It is worth to comment that Theorem~\ref{main} is an extension of what happens
in the classical case $k=1$. There, the same result holds but in statement (b)
only appears one interval. Our second main result will prove that there are
other more significative differences between the case $k=1$ and the cases
$k=2,3.$ These differences are related with the lack of monotonicity of certain
rotation number functions associated to the dynamical systems given by the
Lyness recurrences, see Theorem~\ref{teonomonotonia}. The behaviors of these
rotation number functions are important for the understanding of the
recurrences, because they give the possible periods for them, see
\cite{BR,BC,Z}.

On the other hand in \cite{deA,GKL} it is proved that, at least for some values
of $\{a_n\}_n,$ the behaviour of $\{x_n\}_n$ for the case $k=5$ is totally
different. In particular unbounded positive solutions appear. In the
forthcoming paper \cite{CGM3} we explore in more detail the differences between
the cases $k=1,2,3$ and $k\ge4.$

This paper is organized as follows: Section~\ref{ds} presents the
difference equations that we are studying as discrete dynamical
systems and we state our main results on them, see
Theorems~\ref{rotacions} and \ref{teonomonotonia}. Section~\ref{pt2}
is devoted to the proof of Theorem~\ref{rotacions}. By using it, in
Section~\ref{provaA}, we prove  Theorem~\ref{main} and we give some
examples of how to apply it to determine the number of closed
intervals of the adherence of $\{x_n\}_n$. In
Section~\ref{someproperties} we demonstrate
Theorem~\ref{teonomonotonia} and we also present some examples where
we study in more detail the rotation number function of the
dynamical systems associated to \eqref{eq}.

\section{Main results from the dynamical systems point of view}\label{ds}

In this section we reduce the study of the sequence $\{x_n\}_n$ to the study of
some discrete dynamical systems and we state our main results on them.

First we introduce some notations. When $k=2,$ set
\begin{equation}\label{k=2}a_n\,=\,\left\{\begin{array}{lllr}
a&\mbox{for}&n=2\ell+1,\\ b&\mbox{for}&\,n=2\ell,
\end{array}\right.\end{equation}
and when $k=3,$ set
\begin{equation}\label{k=3}a_n\,=\,\left\{\begin{array}{lllr}
a&\mbox{for}&n=3\ell+1,\\
b&\mbox{for}&n=3\ell+2,\\
c&\mbox{for}&n=3\ell,
\end{array}\right.\end{equation}
where $\ell\in\mathbb{N}$ and $a>0,b>0$ and $c>0.$

We also consider the maps $F_\alpha(x,y)$, with $\alpha\in\{a,b,c\},$ as
$$F_{\alpha}(x,y)=\left(y,\frac{\alpha+y}{x}\right),$$
defined on the open invariant set $\Qm:=\{(x,y):x>0,y>0\}\subset\R^2.$

Consider for instance $k=2.$ The sequence given by (\ref{eq}),
\begin{equation}\label{seqq}
x_1,x_2,x_3,x_4,x_5,x_6,x_7,\ldots,
\end{equation}
can be seen as
\[
(x_1,x_2)\xrightarrow{F_a}(x_2,x_3)\xrightarrow{F_b}(x_3,x_4)
\xrightarrow{F_a}(x_4,x_5)\xrightarrow{F_b}(x_5,x_6)\xrightarrow{F_a}\cdots.
\]
Hence the behavior of (\ref{seqq}) can be obtained from the study of the
dynamical system defined in $\Qm$ by the map:
$$F_{b,a}(x,y):=F_b\circ F_a (x,y)=\left(\frac {a+y}{x},\frac {a+b x+ y}{x
y}\right).
$$
Similarly, for $k=3$ we can consider the map:
$$
F_{c,b,a}(x,y):=F_c \circ F_b\circ F_a(x,y)=\left(\frac {a+bx+y}{xy},\frac
{a+bx+y+cxy}{y \left( a+y \right)} \right).
$$
Notice that both maps have an only fixed point in $\Qm,$ which
depends on $a,b$ (and $c$), that for short we denote  by $\bf{p}.$

It is easy to interpret the invariants for (\ref{eq}) and $k=2,3,$ given in
\cite{JKN,KN}, in terms of first integrals of the above maps, see also
Lemma~\ref{elemental}. We have that
$$V_{b,a}(x,y):=
{\frac {ax^2y+bxy^2+bx^2+ay^2+(b^2+a)x+(b+a^2)y+ab}{xy}},
$$
is a first integral for $F_{b,a}$ and
$$V_{c,b,a}(x,y):= {\frac
{c{x}^{2}y+ax{y}^{2}+b{x}^{2}+b{y}^{2}+(a+bc)x+(c+ab)y+ac}{xy}},$$
is a first integral for $F_{c,b,a}$. The topology of the level sets
of these integrals in $\Qm$ as well as the dynamics of the maps
restricted to them is described by the following result, that will
be proved in Section~\ref{pt2}.

\begin{teo}\label{rotacions}
\begin{itemize}
\item[(i)] The level sets of $V_{b,a}$ (resp. $V_{c,b,a}$) in $Q^+\setminus\{\bf{p}\}$ are
diffeomorphic to circles surrounding $\bf p$, which is the unique fixed point
of $F_{b,a}$ (resp. $F_{c,b,a}$).
\item[(ii)] The action of $F_{b,a}$ (resp. $F_{c,b,a}$) on each level set of
$V_{b,a}$ (resp. $V_{c,b,a}$) contained in $Q^+\setminus\{\bf{p}\}$ is
conjugated to a rotation of the circle.
\end{itemize}
\end{teo}

Once a result like the above one is established the study of the
possible periods of the sequences $\{x_n\}_n$ given by (\ref{eq}) is
quite standard. It suffices, first to get the rotation interval,
which is the open interval formed by all the rotation numbers given
by the above theorem, varying the level sets of the first integrals.
Afterwards, it suffices to find which are the denominators of all
the irreducible rational numbers that belong to the corresponding
interval, see \cite{BC,CGM2,Z}.

The study of the rotation number of these kind of rational maps is
not an easy task, see again \cite{BR,BC,CGM2,Z}. In particular, in
\cite{BR} was proved that the rotation number function parameterized
by the energy levels of the Lyness map $F_a, a\ne1,$ is always
monotonous, solving a conjecture of Zeeman given in \cite{Z}, see
also \cite{ER}. As far as we know, in this paper we give the first
simple example for which this rotation number function is neither
constant nor monotonous. We prove:
\begin{teo}\label{teonomonotonia}
There are positive values of $a$ and $b$, such that the rotation
number function $\rho_{b,a}(h)$ of $F_{b,a}$ associated to the
closed ovals of $\{V_{b,a}=h\}\subset\Qm$ has a local maximum.
\end{teo}

Hence, apart from the known behaviors for the autonomous Lyness maps, that is
global periodicity or monotonicity of the rotation number function, which
trivially holds for $F_{b,a},$ taking for instance $a=b=1$ or $a=b\ne1,$
respectively, there appear more complicated behaviors for the rotation number
function.

Our proof of this result relies on the study of  lower and upper
bounds for the rotation number of $F_{b,a}$ on a given oval of a
level set of $V_{b,a}$ given for some $(a,b)\in(\mathbb{Q}^+)^2$ and
$\{V_{b,a}(x,y)=V_{b,a}(x_0,y_0)\},$ for
$(x_0,y_0)\in(\mathbb{Q}^+)^2$. This can be done because the map on
this oval is conjugated to a rotation and it is possible to use an
algebraic manipulator to follow and to order a finite number
iterates on it, which are also given by points with rational
coordinates. So, only exact arithmetic is used. A similar study
could be done for $F_{c,b,a}$.

\section{Proof of Theorem~\ref{rotacions}}\label{pt2}

{\rec {\it Proof of (i) of Theorem \ref{rotacions}.}} The orbits of
$F_{b,a}$ and $F_{c,b,a}$ lie on the level sets $V_{b,a}=h$ and
$V_{c,b,a}=h$ respectively. These level sets can be seen as the
algebraic curves given by
$$
C_2:=\{c_2(x,y)=ax^2y+bxy^2+bx^2-hxy+ay^2+(b^2+a)x+(b+a^2)y+ab=0\}
$$
and
$$
C_3:=\{c_3(x,y)=cx^2y+axy^2+bx^2-hxy+by^2+(a+bc)x+(c+ab)y+ac=0\},
$$ respectively.

Taking homogeneous coordinates on the projective plane $P\R^2$ both
curves $C_2$ and $C_3$ have the form
$$
C:=\{Sx^2y+Txy^2+Ux^2z+Vxyz+Wy^2z+Lxz^2+Myz^2+Nz^3=0\}.$$

In order to find the branches of them tending to infinity, we examine the
directions of approach to infinity ($z=0$) in the local charts determined by
$x=1$ and $y=1$ respectively.

In the local chart given by $x=1$, the curve $C$ writes as
$$
Sy+Ty^2+Uz+Vyz+Wy^2z+Lz^2+Myz^2+Nz^3=0$$ and it meets the straight line at
infinity $z=0$ when $y(S+Ty)=0.$ Since for both curves $C_2$ and $C_3$ the
coefficients $S$ and $T$ are positive, the only intersection point that could
give points in $\Qm$ is $(y,z)=(0,0).$ The algebraic curve $C$ arrives to
$(y,z)=(0,0)$ tangentially to the line $Sy+Uz=0.$ Since for both curves, $C_2$
and $C_3,$ the coefficients $S$ and $U$ are also positive, we have that the
branches of the level sets tending to infinity are not included in~$Q^+.$

An analogous study can be made in the chart given by $y=1,$ obtaining the same
conclusions.

Moreover, it can be easily checked that in the affine plane both
curves $C_2$ and $C_3$ do not intersect the part of the axes $x=0$
and $y=0$ which is in the boundary of $\Qm$.

In summary, there are no branches of the curves $C_2$ and $C_3$ tending to
infinity or crossing the axes $x=0$ and $y=0$ in $Q^+$, and therefore the
connected components of $C_i\cap Q^+$ for $i=2,3$ are bounded. Notice that this
result in particular already implies the persistence of the sequences given by
\eqref{eq}.

Consider $k=2$. We claim the following facts:

\begin{enumerate}[(a)]
\item In $\Qm$, the set of fixed points of $F_{b,a}$ and the set of singular
points of $C_1$ coincide and they only contain the point ${\bf
p}=(\bar x_,\bar y)$.

\item The function $V_{b,a}(x,y),$ has a
local minimum at $\bf p$.
\end{enumerate}
We remark that item (b) is already known. We present a new simple proof for the
sake of completeness.

From the above claims and the fact that the connected components of the level
sets of $V_{b,a}$ in $\Qm$ are bounded it follows that the level sets of
$V_{b,a}$ in $\Qm\setminus\{\bf{p}\}$ are diffeomorphic to circles.

Let us prove the above claims. The fixed points of $F_{b,a}$ are given by

$$
\begin{cases}
x=\frac{a+y}{x}, \\
y=\frac{a+bx+y}{xy},
\end{cases} \Leftrightarrow
\begin{cases}
x^2=a+y, \\
x(y^2-b)=a+y,
\end{cases}
$$
and so $x^2=x(y^2-b).$ Hence in $Q^+,$ we have that $x=y^2-b$ and the above
system is equivalent to
$$
\begin{cases}
x=y^2-b, \\
xy^2-bx-y-a=0,
\end{cases}
\Leftrightarrow \begin{cases}
x=y^2-b, \\
P(y):=y^4-2by^2-y+b^2-a=0.
\end{cases}
$$
It is not difficult to check that the last system of equations is
precisely the one that gives the critical points of the curves
$V_{a,b}=h.$ Moreover, from the first equation it is necessary that
$x=y^2-b>0$ and hence $y>\sqrt{b}.$ Since $P(y)$ has  only one real
root in $(\sqrt{b},\infty)$ the uniqueness of the critical point
holds.

Let us prove that this critical point corresponds with a local minimum of
$V_{b,a}.$ We will check the usual sufficient conditions given by the Hessian
of $V_{b,a}$ at $\bf p$.

Firstly,
$$
\frac{\partial^2 }{\partial x^2}V_{b,a}(y^2-b,y)=2\,{\frac { \left(
y+a \right) \left( ay+b \right) }{ \left({y}^{ 2}-b \right)
^{3}y}}>0\,\mbox{ for }y>\sqrt{b}.$$ Secondly, the determinant of
the Hessian matrix at the points $(y^2-b,y)$ is
$$
h(y)=\frac{f(y)}{(b-y^2)^4y^4},
$$
where
$$f(y):=(by^2+a-b^2)(-b{y}^{6}+
3\left( a+{b}^{2} \right) {y}^{4}+ 4\left( {a}^{2}+b \right) {y}^{3}+ 3b\left(
2a-{b}^{2} \right) {y}^{2}+b^2({b}^{2}-a ).$$ A tedious computation shows that
$ f(y)=q(y)P(y)+r(y), $ with
\begin{align*} r(y)=&\left(4{a}^{2}{b}^{2}+ 4{a}^{3}+4{b}^{3}+6ab
\right) {y}^{3}+
\left(18{a}^{2}b+ 8{b}^{3}a+3{b}^{2} \right) {y}^{2}\\
&+ \left(-4{b}^{3}{a}^{2}+4{a}^{3}b-4{b}^{4} +12 a{b}^{2}+3{a}^{2} \right)
y-8{b}^{4}a+5{a}^{2}{b}^{2}+3{a}^{3}.
\end{align*}
Observe that if $\bar y$ is the positive root of $P(y)$, then
$\sign(h(\bar y ))=\sign(r(\bar y ))$. Taking into account that
$P(\bar y )=0$ implies that $a=\bar y ^4-2b\bar y ^2-\bar y +b^2$ we
have that
$$r(\bar y )= {\bar y}^{2} \left( 4\,{\bar y}^{3}-4\,b\bar y -1 \right) \left(
b\bar y +1-{\bar y}^{3} \right) ^{2} \left( b-{\bar y}^{2} \right)
^{2}.$$ So, $ \sign\left( 4\,{\bar y}^{3}-4b\bar y
-1\right)=\sign\left(P'(\bar y )\right). $ Since
$P(\sqrt{b})=-\sqrt{b}-a<0,$ $\lim\limits_{y\to\infty} P(y)=+\infty$
and, on this interval, there is only one critical point of $P(y),$
which is simple, we get that $P'(\bar y )>0$ and so $h(\bar y )>0.$
Hence ${\bf p}$ is a local minimum of $V_{b,a}(x,y)$, as we wanted
to prove.

The same kind of arguments work to end the proof for the case $k=3,$ but the
computations are extremely more tedious. We only make some comments.

The fixed points of $F_{c,b,a}$ in $Q^+$ are given by:
\begin{equation*}\label{sistemafcba}
\begin{cases}
P(x):=x^5+ax^4-2x^3-(2a+bc)x^2+(1-b^2-c^2)x+a-bc=0,\, \\
y=Q(x):=\dfrac{(xb+a)}{(x-1)(x+1)}.
\end{cases}
\end{equation*}
It can be proved again that they coincide with the singular points of
$V_{c,b,a}$ in $Q^+.$ This fact follows from the computation of several
suitable resultants between ${\partial V_{c,b,a}}/{\partial x},$ ${\partial
V_{c,b,a}}/{\partial y}$ and $V_{c,b,a}.$

The uniqueness of the fixed point $\bf p$ in $\Qm$ can be shown as follows:
since $Q(x)>0$ implies that $x>1$, we only need to search solutions of
$P_5(x)=0$ in $(1,+\infty)$. With the new variable $z=x-1,$
\begin{equation*}\label{quantessolucions}
\widetilde{P}(z):=P(z+1):=z^5+(c+5)z^4+(8+4c)z^3+(4-ab+4c)z^2-(a+b)z-(a+b)^2=0,
\end{equation*}
Since $\widetilde{P}(0)<0$; $\lim\limits_{z\to +\infty}
\widetilde{P}(z)=+\infty$; and the Descarte's rule, we know that there is only
one positive solution, as we wanted to see.

Finally it can be proved that $\bf p$ is a non-degenerated local
minimum of $V_{c,b,a}$. These computations are complicated, and they
have been performed in a very smart way in \cite{KN}, so we skip
them and we refer the reader to this last reference.\qed

\subsection{Proof of (ii) of Theorem \ref{rotacions}}

In \cite{CGM} it is proved a result that characterizes the dynamics of
integrable diffeomorphisms having a \textsl{Lie Symmetry}, that is a vector
field $X$ such that $ X(F(p))=(DF(p))\,X(p)$. Next theorem states it,
particularized to the case we are interested.

\begin{teo}[\cite{CGM}]\label{teor}
Let $\U\subset \mathbb{R}^2$ be an open set and let $\Phi:\U\rightarrow \U$ be
a diffeomorphism such that:
\begin{enumerate}
\item[(a)] It has a smooth regular first integral $V:\U\rightarrow \R,$ having
its level sets $\Gamma_h:=\{z=(x,y)\in\U\,:\, V(z)=h\}$ as simple
closed curves.
\item[(b)] There exists a smooth function $\mu:\U\rightarrow \R^+$ such
that for any $z\in \U,$
\begin{equation*}\label{mu}
\mu(\Phi(z))=\det(D\Phi(z))\,\mu(z).
\end{equation*}
\end{enumerate}
Then the map $\Phi$ restricted to each $\Gamma_h$ is conjugated to a rotation
with rotation number $\tau(h)/T(h)$, where $T(h)$ is the period of $\Gamma_h$
as a periodic orbit of the planar differential equation
\[
\dot z=\mu(z)\left(-\frac{\partial V(z)}{\partial y},\frac{\partial
V(z)}{\partial x}\right)
\]
and $\tau(h)$ is the time needed by the flow of this equation for
going from any $w\in\Gamma_h$ to $\Phi(w)\in\Gamma_h.$
\end{teo}

Next lemma is one of the key points for finding a Lie symmetry for families of
periodic maps, like the $2$ and $3$--periodic Lyness maps.

\begin{lem}\label{mus}
Let $\{G_a\}_{a\in A}$ be a family of diffeomorphisms of $\U\subset
\mathbb{R}^2$. Suppose that there exists a smooth map $\mu:\U\to \R$
such that for any $a\in A$ and any $z\in \U,$  the equation $
\mu(G_a(z))=\det(DG_a(z))\,\mu(z)$  is satisfied. Then, for every
choice $a_1,\ldots,a_k\in A,$ we have
$$
\mu(G_{[k]}(z))=\det(DG_{[k]}(z))\,\mu(z),$$ where
$G_{[k]}=G_{a_k}\circ\cdots\circ G_{a_{2}}\circ G_{a_1}.$
\end{lem}

\dem It is only necessary to prove the result for $k=2$ because the general
case follows easily by induction. Consider $a_1,a_2\in A$ then
\begin{align*} \mu(G_{a_2,a_1}(z))& =\mu(G_{a_2}\circ
G_{a_1}(z))=\det(DG_{a_2}(G_{a_1}(z)))\,\mu(G_{a_1}(z))= \\
&= \det(DG_{a_2}(G_{a_1}(z)))\,\det(DG_{a_1}(z))\,\mu(z)= \det(D
(G_{a_2}\circ G_{a_1}(z)))\,\mu(z)=\\
&= \det(DG_{a_2,a_1}(z))\mu(z),
\end{align*}
and the lemma follows.\qed

\medskip

{\rec {\it Proof of (ii) of Theorem \ref{rotacions}.}} From part (i) of the
theorem we know that the level sets of $V_{b,a}$ and $V_{c,b,a}$ in
$\Qm\setminus\{{\bf p}\}$ are diffeomorphic to circles. Moreover these
functions are first integrals of $F_{b,a}$ and $F_{c,b,a}$, respectively.
Notice also that for any $a,$ the Lyness map $F_a(x,y)=(y,\frac{a+y}{x})$
satisfies
$$\mu(F_{a}(x,y))=\det(DF_{a}(x,y))\mu(x,y),$$
with $\mu(x,y)=xy.$ Hence, by Lemma \ref{mus},
\[\mu(F_{b,a}(x,y))=\det(DF_{b,a}(x,y))\mu(x,y)\quad\mbox{and}\quad
\mu(F_{c,b,a}(x,y))=\det(DF_{c,b,a}(x,y))\mu(x,y).\] Thus, from
Theorem~\ref{teor}, the result follows.\qed

\medskip

It is worth to comment that once part (i) of the theorem is proved
it is also possible to prove that the dynamics of $F_{b,a}$ (resp.
$F_{c,b,a}$) restricted to the level sets of $V_{b,a}$ (resp.
$V_{c,b,a}$) is conjugated to a rotation by using that they are
given by cubic curves and that the map is birational, see
\cite{JRV}. We prefer our approach because it provides
 a dynamical interpretation of the rotation number together with its
analytic characterization.

\section{Proof of Theorem \ref{main}}\label{provaA}

In order to prove Theorem~\ref{main} we need a preliminary result.
Consider the maps $F_{b,a}$ and $F_{a,b},$ jointly with their
corresponding first integrals $V_{b,a}$ and $V_{a,b}.$ In a similar
way consider $F_{c,b,a}\,,\,F_{a,c,b}$ and $F_{b,a,c}$ with
$V_{c,b,a}\,,\,V_{a,c,b}$ and $V_{b,a,c}.$ Some simple computations
prove the following elementary but useful lemma. Notice that it can
be interpreted as the relation between  the first integrals and the
non-autonomous invariants.

\begin{lem}\label{elemental}
With the above notations:
\begin{itemize}

\item[(i)] $V_{b,a}(x,y)=V_{a,b}(F_a(x,y)).$

\item[(ii)] $V_{c,b,a}(x,y)=V_{a,c,b}(F_a(x,y))=V_{b,a,c}(F_b(F_a(x,y))).$
\end{itemize}
\end{lem}

{\rec {\it Proof of Theorem \ref{main}.}} We split the proof in two
steps. For $k=2,3$ we first   prove that there are only two types of
behaviors for $\{x_n\}_n$, either this set of points is formed by
$kp$ points for some positive integer $p,$ or it has infinitely many
points whose adherence is given by at most $k$ intervals. Secondly,
in this later case, we provide an algebraic way for studying the
actual number of intervals.

\noindent {\bf First step:} We start with the case $k=2.$ With the notation
introduced in \eqref{k=2}, it holds that
\[
F_{b,a}(x_{2n-1},x_{2n})=(x_{2n+1},x_{2n+2}), \quad
F_{a,b}(x_{2n},x_{2n+1})=(x_{2n+2},x_{2n+3}),
\]
where  $(x_1,x_2)\in\Qm$ and $n\ge 1$.  So the  odd terms of the sequence
$\{x_n\}_n$  are contained in the projection on the $x$-axis of the oval of
$\{V_{b,a}(x,y)=V_{b,a}(x_1,x_2)=h\}$ and the even terms  in the corresponding
projection of $\{V_{a,b}(x,y)=V_{a,b}(F_{a}(x_1,x_2))=h\}$, where notice that
we have used Lemma~\ref{elemental}.

Recall that the ovals of $V_{b,a}$  are invariant by $F_{b,a}$ and
the ovals of $V_{a,b}$  are invariant by $F_{a,b}$. Notice also that
the trivial equality $F_a\circ F_{b,a}=F_{a,b}\circ F_a$ implies
that the action of $F_{b,a}$ on $\{V_{b,a}(x,y)=h\}$ is conjugated
to the action of $F_{a,b}$ on $\{V_{a,b}(x,y)=h\}$ via $F_a.$

From Theorem \ref{rotacions} we know that $F_{b,a}$ on the
corresponding oval is conjugated to a rotation of the circle. Hence,
if the corresponding rotation number is rational, then the orbit
starting at $(x_1,x_2)$ is periodic, of period say $q,$ then  the
sequence $\{x_n\}_n$ is $2q$-periodic. On the other hand  if the
rotation number is irrational, then the orbit of $(x_1,x_2)$
generated by $F_{b,a}$ fulfills densely the oval of
$\{V_{b,a}(x,y)=h\}$ in $\Qm$ and hence the subsequence of odd terms
also fulfills densely the projection of $\{V_{b,a}(x,y)=h\}$ in the
$x$-axis.  Clearly, the sequence  of even terms do the same with the
projection of the oval of $\{V_{a,b}(x,y)=h\}.$

Similarly when $k=3$ the equalities
\begin{align*}
F_{c,b,a}(x_{3n-2},x_{3n-1})&=(x_{3n+1},x_{3n+2}), \\
F_{a,c,b}(x_{3n-1},x_{3n})&=(x_{3n+2},x_{3n+3}),\\
F_{b,a,c}(x_{3n},x_{3n+1})&=(x_{3n+3},x_{3n+4}),
\end{align*}
where $n\ge1$, allow to conclude that each term $x_m,$ of the sequence
$\{x_n\}_n$ where we use the notation (\ref{k=3}), is contained in one of the
projections on the $x$-axis of the ovals $\{V_{c,b,a}(x,y)=
V_{c,b,a}(x_1,x_2)=:h\}$ and $\{V_{a,c,b}(x,y)=h\}$ and $\{V_{b,a,c}(x,y)=h\}$,
according with the remainder of $m$ after  dividing it by 3. The rest of the
proof in this case follows as in the case $k=2.$ So the first step is done.

\noindent {\bf Second step:} From the above results it is clear that
the problem of knowing the number of connected components of the
adherence of $\{x_n\}_n$ is equivalent to the control of  the
projections of several invariant ovals on the $x$-axis. The strategy
for $k=3$, and analogously for $k=2$, is the following. Consider the
ovals contained in the level sets given by $\{V_{c,b,a}(x,y)=h\}$,
$\{V_{a,c,b}(x,y)=h\}$ and $\{V_{b,a,c}(x,y)=h\}$ and denote by
$I=I(a,b,c,h), J=J(a,b,c,h)$ and $K=K(a,b,c,h)$ the corresponding
closed intervals of $(0,\infty)$ given by their projections on the
$x$-axis.

We want to detect the values of $h$ for which two of the intervals,
among $I,J$ and $K,$ have exactly one common point. First  we seek
their boundaries. Since the level sets are given by cubic curves,
that are quadratic with respect the $y$-variable, these points will
correspond with values of $x$ for which the discriminant of the
quadratic equation with respect to $y$ is zero. So, we compute
\begin{align*}
R_1(x,h,a,b,c)&:=\mbox{dis}\,(xyV_{c,b,a}(x,y)-hxy,y)=0,\\
R_2(x,h,a,b,c)&:=\mbox{dis}\,(xyV_{a,c,b}(x,y)-hxy,y)=0,\\
R_3(x,h,a,b,c)&:=\mbox{dis}\,(xyV_{b,a,c}(x,y)-hxy,y)=0.
\end{align*}
Now we have to search for relations among $a,b,c$ and $h$ for which
 two of these three functions have some common solution, $x.$ These
 relations can be obtained by computing some suitable resultants.

Taking the resultants of $R_1$ and $R_2$; $R_2$ and $R_3$; and $R_1$
and $R_3$ with respect $x$ we obtain three polynomial equations
$R_4(h,a,b,c)=0$, $R_5(h,a,b,c)=0$ and $R_6(h,a,b,c)=0.$ In short,
once $a,b$ and $c$ are fixed we have obtained three polynomials in
$h$ such that  a subset of their zeroes give the bifurcation values
which separate the number of intervals of the adherence of
$\{x_n\}_n$. See the results of Proposition~\ref{ex2} and
Example~\ref{ex3} for concrete applications of  the method.

Before ending the proof we want to comment that for most values of
$a$, $b$ and $c,$ varying $h$ there appear the three possibilities,
namely $1$, $2$ or $3$ different intervals. The last case appears
for values of $h$ near $h_c:=V_{c,b,a}({\bf p})$,  because the first
 coordinates of the three points $\bf p,$ $F_a({\bf p})$
 and $F_b(F_a({\bf p}))$ almost never coincide. The other situations
can be obtained  by increasing $h.$ \qed

\begin{propo}\label{ex2} Consider the recurrence \eqref{eq} with $k=2$ and
 $\{a_n\}_n$ as in \eqref{k=2}
taking  the values $a=3$ and $b=1/2.$ Define
$h_c={(12z^3-33z+7)}/{(2(z^2-3))}\simeq 17.0394, $ where $z\simeq
2.1513$ is the biggest positive real root of \, $2z^4-12z^2-2z+17,$
and $h^*\simeq 17.1198,$ is the smallest positive root of
\[p_4(h):=112900h^4-2548088h^3-48390204h^2+564028596h+7613699255.\]
 Then,
\begin{enumerate}[(i)]

\item The initial condition $(x_1,x_2)=(z,z^2-3)$ gives a two periodic
recurrence $\{x_n\}_n$. Moreover $V_{1/2,3}(z,3-z^2)=h_c.$

\item Let $(x_1,x_2)$ be any positive initial conditions, different
from $(z,z^2-3),$ and set $h=V_{1/2,3}(x_1,x_2).$ Let $\rho(h)$
denote the rotation number of $F_{1/2,3}$ restricted to the oval of
$\{V_{1/2,3}(x,y)=h\}.$ Then
\begin{enumerate}[(I)]

\item If $\rho(h)=p /q\in\mathbb{Q},$ with $\gcd(p,q)=1,$ then the sequence
$\{x_n\}_n$ is $2q$-periodic.

\item If $\rho(h)\not\in\mathbb{Q}$ and $h\in (h_c,h^*)$ then the adherence of
the sequence $\{x_n\}_n$ is formed by two disjoint closed intervals.

\item If $\rho(h)\not\in\mathbb{Q}$ and $h\in [h^*,\infty)$ then the adherence of
the sequence $\{x_n\}_n$ is one closed interval.

\end{enumerate}

\end{enumerate}
\end{propo}

We want to remark that, from a computational point of view, the case
(I) almost never is detected. Indeed, taking  $a$ and $b$ rational
numbers and  starting with rational initial conditions, by using
Mazur's theorem it can be seen that almost never the rotation number
will be rational, see the proof of  \cite[Prop. 1]{BR}. Therefore,
in numeric simulations only situations (II) and (III) appear, and
the value $h=h^*$ gives the boundary between them. In general, for
$k=2,$ the value $h^*$ is always the root of a polynomial of degree
four, which is constructed from the values of $a$ and $b$.

\noindent{\it Proof of Proposition~\ref{ex2}.}   Clearly $(z,3-z^2)$
is the fixed point of $F_{b,a}$ in $\Qm.$ Some computations give the
compact expression of $h_c:=V_{a,b}(z,z^2-3).$ To obtain the values
$h^*$ we proceed as in the proof of Theorem~\ref{main}. In general,
\begin{align*}
R_1(x,h,a,b)&:=\mbox{dis}\,(xyV_{b,a}(x,y)-hxy,y)\\
&\phantom{:}=(ax^2-hx+a^2+b)^2-4(bx+a)(bx^2+b^2x+ax+ab),\\
R_2(x,h,a,b)&:=\mbox{dis}\,(xyV_{a,b}(x,y)-hxy,y)\\
&\phantom{:}=(bx^2-hx+a+b^2)^2-4(ax+b)(ax^2+a^2x+bx+ab).
\end{align*}
Then we have to compute the resultant of the above polynomials with
respect to $x$. It always decomposes as the product of two quartic
polynomials in $h.$ Its expression is very large, so we only give it
when $a=3$ and $b=1/2$. It writes as
\[
\frac{625}{65536}\left(4h^4-1176h^3+308h^2+287380h+1816975\right)p_4(h).
\]
It has four real roots, two for each polynomial. Some further work
proves that the one that interests us is the smallest one of $p_4.$
\qed

We also give an example when $k=3,$ but skipping all the details.

\begin{example}\label{ex3} Consider the recurrence \eqref{eq} with $k=3$ and
 $\{a_n\}_n$ as in \eqref{k=3}
taking  the values $a=1/2, b=2$ and $c=3.$ Then  for any positive
initial conditions $x_1$ and $x_2$,  $V_{c,b,a}(x_1,x_2)=h\ge
V_{c,b,a}({\bf p})=h_c\simeq 15.9283.$ Moreover if the rotation
number of $F_{c,b,a}$ associated to the oval $\{V_{c,b,a}(x,y)=h\}$
is irrational then the adherence of $\{x_n\}_n$ is given by:
\begin{itemize}
\item Three intervals when $h\in (h_c,h^*),$ where $h^*\simeq
15.9614;$
\item Two intervals when $h\in [h^*,h^{**}),$ where $h^{**}\simeq
16.0015;$
\item One interval when $h\in[h^{**},\infty).$
\end{itemize}
The values $h^*$ and $h^{**}$ are roots of two polynomials of degree
8 with integer coefficients that can be explicitly given.
\end{example}

\section{Some properties of the rotation number function
}\label{someproperties}

From Theorem~\ref{rotacions} it is natural to introduce the {\it
rotation number function} for $F_{b,a}$:
\[
\rho_{b,a}:[h_c,\infty)\longrightarrow (0,1),
\]
where $h_c:=V_{b,a}({\bf p})$, as the map that associates to each
invariant oval $\{V_{b,a}(x,y)=h\}$, the rotation number
$\rho_{b,a}(h)$ of the function $F_{b,a}$ restricted to it. The
following properties hold:
\begin{enumerate}[(i)]
\item The function $\rho_{b,a}(h)$ is analytic for $a>0,b>0$,
$h>h_c$ and it is continuous at $h=h_c.$ This can be proved from the
tools introduced in  \cite[Sec. 4]{CGM2}.

\item The value $\rho_{b ,a}(h_c)$ is given by the argument over $2\pi$ of the
eigenvalues (which have modulus one due to the integrability of
$F_{b,a}$) of the differential of $F_{b,a}$ at ${\bf p}$.

\item $\rho_{b,a}(h)=\rho_{a,b}(h).$

\item $\rho_{a,a}(h)=2\rho_a(h) \mod 1,$ where $\rho_a$ is the rotation
number\footnote{Notice that given a map of the circle there is an
ambiguity between $\rho$ and $1-\rho$ when one considers its
rotation number. So, while for us the rotation number of the
classical Lyness map for $a=1$ is $4/5$, in other papers it is
computed as $1/5.$}  function associated to the classical Lyness
map. Then, from the results of \cite{BC} we know that
$\rho_{1,1}(h)\equiv3/5,$  that for $a\ne1,$ positive,
$\rho_{a,a}(h)$ is monotonous and
$\lim_{h\to\infty}\rho_{a,a}(h)=3/5.$
\end{enumerate}

Note that item (iii) follows because $F_{a,b}$ is conjugated with
$F_{b,a}$ via $\psi=F_a$ which is a diffeomorphism of $\Qm$, because
$ \psi^{-1} F_{a,b} \psi=F_a^{-1}F_a F_b F_a=F_b F_a=F_{b,a}. $
Since $\psi$ preserves the orientation, the rotation number
functions of $F_{a,b}$ and $F_{b,a}$ restricted to the corresponding
ovals must coincide.

Similar results to the ones given above hold for $F_{c,b,a}$ and its
corresponding rotation number function. Some obvious differences
are:
\begin{align*}
&\rho_{c,b,a}(h)=\rho_{b,a,c}(h)=\rho_{a,c,b}(h)\,
&&\rho_{a,a,a}(h)=3\rho(h)\,\mod1,\\
&\rho_{1,1,1}(h)=2/5,  &&\lim_{h\to\infty}\rho_{a,a,a}(h)=2/5.
\end{align*}

We are convinced that when $a>0$ and $b>0,$
\[\lim_{h\to\infty}\rho_{b,a}(h)=3/5\quad\mbox{and}\quad\lim_{h\to\infty}\rho_{c,b,a}(h)=2/5,
\]
but we have not been able to prove these equalities. If they were
true, by combining them with the values of the rotation number
function at $h=h_c$ this would give a very useful information to
decide if, apart from the trivial cases $a=b=1 (c=1),$ there are
other cases for which the rotation number function is constant.
Notice that in these situations  the maps $F_{b,a}$ or $F_{c,b,a}$
would be globally periodic in $\Qm.$ This information, together with
the values at $h_c$, also would be useful to know the regions where
the corresponding functions could be increasing or decreasing.
Finally notice that this rotation number at infinity is not
continuous when we approach to $a=0$ or $b=0$, where the recurrence
and the first integral are also well defined on $\Qm.$ For instance
$\rho_{0,0}(\rho)\equiv 2/3$ and the numerical experiments of next
subsection seem to indicate that for $a>0$ or $b>0,$
\[
 \lim_{h\to\infty}\rho_{0,a}(h)=\lim_{h\to\infty}\rho_{b,0}(h)=5/8.
\]

Before proving Theorem~\ref{teonomonotonia} we introduce with an
example the algorithm that we will use along this section to compute
lower and upper bounds for the rotation number. We have implemented
it in an algebraic manipulator. Notice also that when we apply it
taking rational values of $a$ and $b$ and rational initial
conditions, it can be used as a method to achieve proofs, see next
example or the proof  of Theorem~\ref{teonomonotonia}.

\begin{figure}[h]
\begin{center}
\includegraphics[scale=0.35]{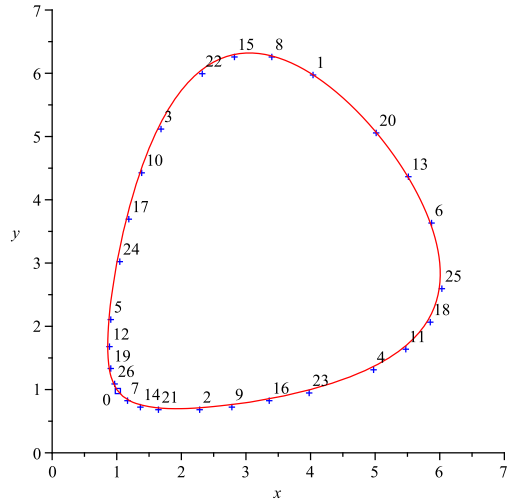}
\end{center}
\caption{Oval of $\{V_{2,3}(x,y)=34\}$ with 27 iterates of
$F_{2,3}$. The label $0$ indicates the initial condition $(1,1)$,
and the label $k, k=1,\ldots,26,$ corresponds with the $k$-th point
of the orbit.}\label{fig1}
\end{figure}

Fix $a=3,$ $b=2$ and $(x_0,y_0)=(1,1).$ Then $h=V_{2,3}(1,1)=34.$
Compute for instance the 27 points of the orbit starting at $(1,1),$
\[
(x_1,y_1)=(4,6),\quad (x_2,y_2)=\left(\frac9 4,
\frac{17}{24}\right),\quad (x_3,y_3)=\left(\frac{89}{54},
\frac{788}{153}\right),\ldots
\]
and consider them as  points on the oval $\{V_{2,3}(x,y)=34\},$ see
Figure~\ref{fig1}.

\begin{figure}[h]
\begin{center}
\includegraphics[height=5 cm,width=12 cm]{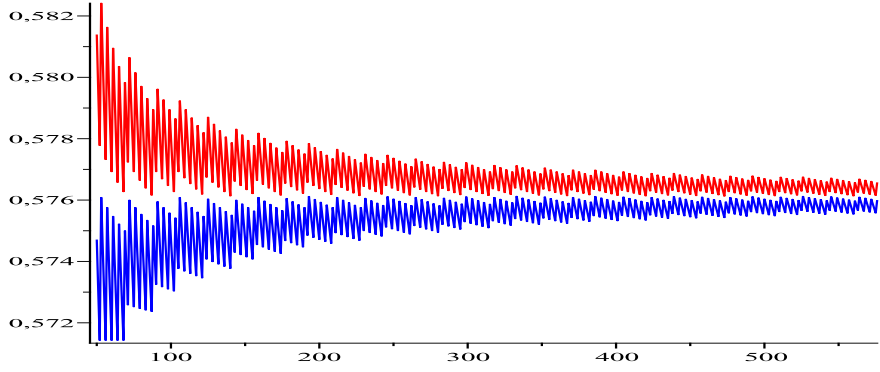}
\end{center}
\caption{Lower and upper bounds for $\rho_{2,3}(34)$ obtained after
following  some points of the orbit starting at
$(1,1)$.}\label{fig2}
\end{figure}

We already know that the restriction of $F_{2,3}$ to the given oval
is conjugated to a rotation, with rotation number
$\rho:=\rho_{2,3}(34)$ that we want to estimate. This can be done by
counting the number of turns that give the points $(x_j,y_j)$, after
fixing some orientation in the closed curve. We orientate the curve
with the counterclockwise sense. So, for instance we know that the
second point $(x_2,y_2)$ has given more that one turn and less than
two, giving that $1<2\rho<2,$ and hence that $\rho\in(1/2,1).$ Doing
the same reasoning with all the points computed we obtain,
\begin{align*}
&4<7\rho<5&\Rightarrow&&& \rho\in\left(\frac47,\frac57\right),\\
&8<14\rho<9&\Rightarrow&&& \rho\in\left(\frac8{14},\frac9{14}\right),\\
&10<19\rho<11&\Rightarrow&&& \rho\in\left(\frac{10}{19},\frac{11}{19}\right),\\
 &14<26\rho<15&\Rightarrow&&&
\rho\in\left(\frac{14}{26},\frac{15}{26}\right),
\end{align*}
where we have only written the more relevant informacions obtained,
which are given by the points of the orbit closer to the initial
condition. So, we have shown that
\[
0.5714\simeq\frac{4}{7}<\rho_{2,3}(34)<\frac{15}{26}\simeq 0.5769.
\]

In Figure~\ref{fig2} we represent several successive lower and upper
approximations obtained while the orbit is turning around the oval.
We plot around six hundred steps, after skipping the first fifty
ones. By taking 1000 points we get
\[
0.5761246\simeq\frac{338}{578}<\rho_{2,3}(34)<\frac{473}{821}\simeq
0.5761267,
\]
and after 3000 points,
\[
0.57612457\simeq\frac{338}{578}<\rho_{2,3}(34)<\frac{1472}{2555}\simeq
0.57612524.
\]
In fact when we say that
$\rho_{2,3}(34)\in(\rho_{\mathrm{low}},\rho_{\mathrm{upp}})$, the
value $\rho_{\mathrm{low}}$ is the upper lower bound obtained by
following all the considered points  of the orbit, and
$\rho_{\mathrm{upp}}$ is the lowest upper bound. Notice that taking
1000 or 3000 points we have obtained the same lower bound for
$\rho_{2,3}(34).$

Let us prove  Theorem~\ref{teonomonotonia} by using the above
approach.

\bigskip

\noindent{\it Proof of Theorem \ref{teonomonotonia}.} Consider
$a=1/2,$  $b=3/2$ and the three points
\[
\mathbf{p}^1=\left(\displaystyle{\frac{149}{100}},\displaystyle{\frac{173}{100}}
\right),\quad\mathbf{p}^2=\left(\displaystyle{\frac{3}{40}},\displaystyle{\frac{173}{100}}
\right),\quad\mathbf{p}^3=\left(\displaystyle{\frac{1}{1000}},\displaystyle{\frac{173}{100}}
\right).
\]
Notice that
\begin{align*}
&h_1:=V_{3/2,1/2}(\mathbf{p}^1)=\frac{10655559}{1288850}\simeq8.27,\\
&h_2:=V_{3/2,1/2}(\mathbf{p}^2)=\frac{9328327}{207600}\simeq44.93,\\
&h_3:=V_{3/2,1/2}(\mathbf{p}^3)=\frac{1056238343}{346000}\simeq3052.71.
\end{align*}
Hence $h_c<h_1<h_2<h_3.$ By applying the algorithm described above,
using 100 points of each orbit starting at each ${\bf p}^j,
j=1,2,3,$ we obtain that
\[
\rho_{3/2,1/2}(h_1),\rho_{3/2,1/2}(h_3)\in\left(\frac
35,\frac{59}{98}\right)\quad\mbox{and}\quad
\rho_{3/2,1/2}(h_2)\in\left(\frac{56}{93},\frac{53}{88}\right).
\]
Since $59/98<56/93$ we have proved that the function
$\rho_{3/2,1/2}(h)$ has at least a local maximum in $(h_1,h_3).$
From the continuity of the rotation number function, with respect
$a,b$ and $h$, we notice that this result also holds for all values
of $a$ and $b$ in a neighborhood of $ a=1/2,b=3/2.$ \qed

\bigskip

We believe that with the same method it can be proved that a similar
result to the one given in  Theorem~\ref{teonomonotonia} holds for
some maps $F_{c,b,a},$ but we have decided do not perform this
study.

\subsection{Some  numerical explorations for $k=2.$}\label{seccionumeric}

\begin{figure}[h]
\begin{center}
\includegraphics[height=7 cm,width=14 cm]{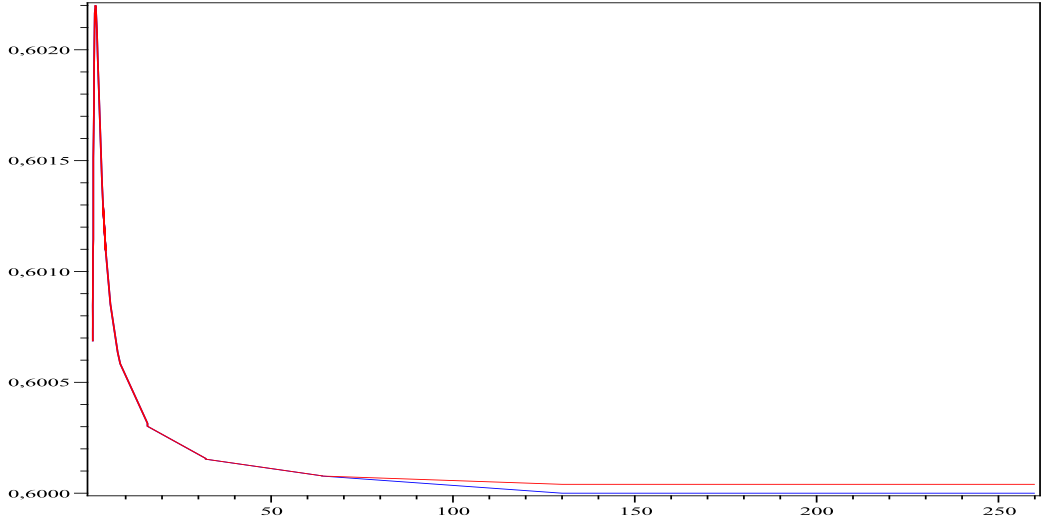}
\end{center}
\caption{Lower and upper bounds for $\rho_{3/2,1/2}(h)$. On the
horizontal axis we represent $-\log_{10}(h)$ and, on the vertical
axis, the value of the rotation number. Notice that  for values of
$-\log_{10}(h)$ smaller that 70 both values are indistinguishable in
the Figure.}\label{fig3}
\end{figure}

We start by studying with more detail the rotation number function
$\rho_{3/2,1/2}(h)$, that we have considered to prove
Theorem~\ref{teonomonotonia}. In this case the fixed point is ${\bf
p}\simeq(1.493363282,1.730133891)$ and
$h_c=V_{b,a}({\mathbf{p}})=8.267483381.$ 
Moreover $\rho_{b,a}({h}_c)\simeq0.6006847931$. By applying our
algorithm for approximating the rotation number, with $5000$ points
on each orbit, we obtain the results presented in Table~1. In
Figure~\ref{fig3} we also plot the upper and lower bounds of
$\rho_{3/2,1/2}(h)$ that we have obtained by using a wide range of
values of $h.$

\begin{center}
\vglue 0.2cm

\begin{tabular}{|r|r|r|r|}
\hline Init. cond. $(x,\bar y)$ & Energy level $h$ &
$\rho_{\mathrm{low}}(h)\qquad$ & $\rho_{\mathrm{upp}}(h)
\qquad$\\
\hline

$\bar x$&$h_c\simeq 8.2675$& $\simeq 0.6006848$ & $\simeq 0.6006848$ \\

1.3&$8.3068$& $\frac{173}{288}\simeq 0.6006944$ & $\frac{2938}{4891}\simeq 0.6006951$ \\

0.75&$9.2747$& $\frac{1435}{2388}\simeq 0.6009213$ & $\frac{2087}{3473}\simeq 0.6009214$ \\

0.3&$14.7566$& $\frac{1548}{2573}\simeq 0.6016323$ & $\frac{2285}{3798}\simeq 0.6016324$ \\

0.075&$44.9347$& $\frac{657}{1091}\simeq 0.6021998$ & $\frac{2354}{3909}\simeq 0.6022001$ \\

0.001&$3052.75$& $\frac{2927}{4867}\simeq 0.6013972$ & $\frac{86}{143}\simeq 0.6013986$ \\

$5\cdot 10^{-6}$&$609716.07$& $\frac{1832}{3049}\simeq 0.6008527$ & $\frac{1409}{2345}\simeq 0.6008529$ \\

$5\cdot 10^{-256}$&$6.097\cdot 10^{255}$& $\frac{3}{5}= 0.6$ & $\frac{2999}{4998}\simeq 0.6000400$ \\
\hline
\end{tabular}
\vglue 1 cm

Table 1: Lower and upper bounds of the rotation number
$\rho_{3/2,1/2}(h)$, for some orbits of $F_{3/2,1/2}$ starting at
$(x,\bar y),$ where ${\bf p}=(\bar x,\bar y).$

\end{center}

For other values of $a$ and $b$ we obtain different behaviors. All
the experiments are performed by starting at the fixed point ${\bf
p}=(\bar x,\bar y),$ and increasing the energy level by taking
initial conditions of the form $(x,\overline{y})$, by decreasing $x$
to $0$. With this process we take orbits approaching to the boundary
of $\Qm$, that is lying on level sets of $V_{b,a}$ with increasing
energy. The step in the decrease of $x$ (and therefore in the
increase of $h$) is not uniform, and it has been manually tuned
making it smaller in those regions where a possible non monotonous
behavior could appear.

Consider the set of parameters $\Gamma=\{(a,b),\in [0,\infty)^2\}$,
where notice that we also consider the boundaries $a=0$ or $b=0,$
where   the map $F_{b,a}$ is   well defined. We already know that
the rotation number function behaves equal at $(a,b)$ and $(b,a).$
Moreover we know perfectly its behavior on the diagonal $(a,a)$
(when $a<1$ it is monotonous decreasing and when $a>1$ it is
monotonous increasing) and that $\rho_{1,1}(h)\equiv 4/5$ and
$\rho_{0,0}(h)\equiv 2/3.$ Hence a good strategy for an numerical
exploration can be to produce sequences of experiments using our
algorithm by fixing some $a\ge0$ and varying $b$. For instance we
obtain:

\begin{itemize}

\item Case $a=1/2.$ For all the values of $b>0$ considered, the rotation number function
seems to tend to $3/5$ when $h$ goes to infinity. Moreover it seems

\begin{itemize}

\item monotone decreasing for  $b\in\{1/4,1\};$

\item  to have a unique maximum when $b\in\{7/5,3/2\}$;

\item monotone increasing for $b\in\{2,3\}.$

\end{itemize}

\item Case $a=0.$ For all the values of $b>0$ considered, the rotation
number function seems to tend to $5/8$ when $h$ goes to infinity.
Moreover it seems

\begin{itemize}

\item monotone decreasing for    $b\in\{1/10,3/10,1/2\}$;

\item  to have a unique maximum when   $b\in\{7/10,3/4\}$;

\item  monotone decreasing for  $b\in\{1,5\}$.

\end{itemize}

\end{itemize}

The  above results, together with some other experiments for other
values of $a$ and $b$, not detailed in this paper,  indicate the
existence of a subset of positive measure in  $\Gamma$ where the
corresponding rotation number functions seem to present an unique
maximum. This subset probably separates two other subsets
of~$\Gamma$, one where $\rho_{b,a}(h)$ is monotonically decreasing
to $3/5$ , and another one where $\rho_{b,a}(h)$ increases
monotonically to the same value. The ``oscillatory subset'' seems to
shrink to $(a,b)=(1,1)$ when it approaches to the line $a=b$ and
seems to finish in one interval on each of the borders $\{a=0\}$ and
$\{b=0\}$. Further analysis must be done in this direction in order
to have a more accurate knowledge of  the bifurcation diagram
associated to the behavior of  $\rho_{b,a}$  on $\Gamma$.

\end{document}